\documentclass[twoside,12pt,leqno]{article}
\advance\topmargin by -\headheight\advance\topmargin by -\headsep

\usepackage{amsmath}
\usepackage{amsthm}
\usepackage{amssymb}
\usepackage{amscd}
\usepackage{graphicx}
\usepackage{epsfig}

\numberwithin{equation}{section} \allowdisplaybreaks
\newtheorem{theorem}{Theorem}[section]
\newtheorem{lemma}{Lemma}[section]
\newtheorem{proposition}{Proposition}[section]
\newtheorem{corollary}{Corollary}[section]

\theoremstyle{definition}
\newtheorem{definition}{Definition}[section]

\newtheorem{remark}{Remark}[section]

\begin{document}



\setcounter{page}{1}
\thispagestyle{empty}



\markboth{Y. Kosmann-Schwarzbach and C. Laurent-Gengoux}{The modular
  class of a twisted Poisson structure}

\label{firstpage}

\bigskip

\bigskip

\centerline{{\Large The modular class of a twisted Poisson structure}}

\bigskip
\bigskip
\centerline{{\large by Yvette Kosmann-Schwarzbach and Camille Laurent-Gengoux}}

\vspace*{.7cm}

\begin{abstract}
We study the geometric and algebraic properties of the twisted Poisson
structures on Lie algebroids,
leading to a definition of their modular class and to
an explicit determination of a representative
of the modular class, in particular
in the case of a twisted Poisson manifold.\medskip
\\
{\it R\'esum\'e.}
\noindent Nous \'etudions des propri\'et\'es g\'eom\'etriques et alg\'ebriques
des structures de Poisson tordues sur les alg\'ebro{\"{\i}}des de Lie,
permettant une d\'efinition de leur classe modulaire et une d\'etermination
explicite d'un repr\'esentant de celle-ci, en particulier
pour le cas d'une vari\'et\'e de
Poisson tordue.
\end{abstract}

\pagestyle{myheadings}

\section{Introduction}

The primary aim of this paper is to extend the definition and
properties of the modular class of Poisson manifolds to the case of
manifolds with a {\it twisted Poisson structure}. Moreover
we show that the notion of modular class
can be extended to the case of
Lie algebroids with a twisted Poisson structure.

There are many ways in which the Jacobi identity for a skew-symmetric
bracket can be violated. On a manifold with a
twisted Poisson structure, the Jacobi identity for the
Poisson bracket holds only up to an additional term
involving a closed 3-form, called the
{\it background 3-form}.
Such structures appeared in string theory.
The quantization of an open membrane coupled to a background 3-form
was interpreted in the work of
Jae-Suk Park \cite{P}
as a deformation of the theory
without background;
the defining condition, equation (\ref{3forme}) below,
appears as the quantum master equation in the Batalin-Vilkovisky
quantization of an action functional containing a term, the $C$-field,
which  is a closed 3-form.
Meanwhile, in \cite{CS}, Cornalba and Schiappa introduced a
star-product deformation which is not only non-commutative but also
non-associative; the Jacobi identity
for the associated commutator bracket was
therefore violated, and that identity took the form that appears in formula
\eqref{eq:Jacobi} below.
For Klim\v cik and Strobl \cite{strobl},
equation (\ref{3forme}) appears as the condition for the constraints
for a Lagrangian system to be first-class, i.e., to span a subalgebra
of the Poisson algebra,
when the action is that of a Poisson $\sigma$-model
to which is added a term analogous to the term of Wess and Zumino as,
e.g., in the Wess-Zumino-Witten model. For this reason,
Klim\v cik et Strobl proposed to call such structures
{\it WZW-Poisson structures} or {\it WZ-Poisson structures}.
``Poisson geometry with a 3-form background'' was then
studied by \v Severa and Weinstein \cite{severa}
who showed that such a structure is
a Dirac structure in a Courant algebroid whose bracket is
defined by means of the 3-form. They called it
a {\it Poisson structure with background} or
{\it a twisted Poisson structure}, the term that we have adopted here, in
spite of a possible confusion of terminology
evoked at the very end of our paper.
While
any twisted  Poisson manifold is locally equivalent to a
genuine Poisson manifold, global phenomena make this generalized case
interesting.

The {\it modular vector fields}
of Poisson manifolds already figured in Koszul's
1985 article
\cite{Ko}, and some of their properties and applications appear in the
work of Dufour
and Haraki \cite{DH}, who called them ``curl'' (rotationnel, in
French), and in other papers of the early nineties.
Weinstein, in \cite{W}, related this notion to the modular automorphism
group of von Neumann algebras, gave it the name that has been adopted
in the literature, and introduced the notion of modular class.

Given an orientable Poisson manifold, choose a volume form and
associate to each smooth function the divergence of its hamiltonian
vector field. The map thus obtained is a derivation and is, by
definition, the modular vector field.
The basic observation is that this vector field is closed in the Poisson
cohomology and that its cohomology class, the {\it modular class},
does not depend on the choice of the
volume form. The non-orientable case can be dealt with, replacing
volume forms by densities. Further advances, already announced in
\cite{W},
were made in the article
of Evens, Lu and Weinstein \cite{ELW} where
the modular class for a Lie algebroid was defined and where it was
shown that the modular class of a Poisson manifold was one half that of its
cotangent Lie algebroid.
At the same time, Huebschmann developed a powerful algebraic theory
in the framework of Lie-Rinehart algebras \cite{H} \cite{H2}
 which recovered the results
of \cite{ELW} when applied to the case where a Lie-Rinehart algebra
is the space of sections of a Lie algebroid.
Duality properties were proved by these authors and by Xu in \cite{Xu}.

In this article, we shall follow the approach of \cite{yks2000},
 where the
modular vector fields are characterized in terms of the difference of
 two generating operators of square zero of the Gerstenhaber algebra
 associated to the given Lie algebroid.
After brief preliminary results on operators of order 1 and 2 on graded
 algebras,
we introduce, in Section \ref{section3},
operators on forms and multivectors and a vector field,
defined in terms of a bivector and a 3-form, that are the needed
 ingredients of the construction of the modular class.
This makes sense on an
 arbitrary vector bundle.
In Section~\ref{section4}, we define the Lie algebroids with a twisted
 Poisson structure, which comprise the cotangent Lie algebroids of
 twisted Poisson manifolds and the triangular Lie bialgebroids.
We show in Section \ref{section5} that on
a Lie algebroid with a twisted Poisson structure, there
 exist generators of square zero of the Gerstenhaber algebra of
multivectors, defined in terms of the operators of Section
 \ref{section3},
and the definition and properties of the modular vector field
(Section~ \ref{sectionZ}) follow. That vector field
is closed in the Lie algebroid cohomology
 and its class is
 well-defined (Theorem \ref{theorem6.1}). It is the sum of the vector field
 $X_{\pi,\lambda}$, depending on the bivector~$\pi$ and the volume
 form $\lambda$ (or density in the non-orientable case), that
 appears in the untwisted case, and which is no longer closed in the
 twisted case, and the vector field, $Y_{\pi,\psi}$,
depending on the bivector $\pi$
and the 3-form $\psi$ defining the twisted Poisson structure.
In Section \ref{sectionELW},
we then recall the construction of Evens, Lu and Weinstein \cite{ELW}
and prove that, as expected,  when the Lie algebroid is the
 tangent bundle of a twisted Poisson manifold, $M$, the class that
 we have defined is one half
 the class, as defined in
 \cite{ELW}, of the cotangent Lie algebroid $T^*M$.
The examples described in Section \ref{section8}
show that this is not
the case in general, even for Lie algebras considered as Lie algebroids over
 a point. As another example, we show that the modular classes 
of the Lie groups equipped with
the twisted Poisson structure
 introduced in  \cite{severa} vanish.
Many of the features of the usual Poisson case can be recovered but
 new phenomena appear in the case of the Lie algebroids with a twisted
 Poisson structure.

A full understanding of the relationship between the modular class
that we introduce and the modular classes defined
in \cite{W} \cite{ELW} for Lie algebroids, and further
justification for the generalization that we propose are
provided by the consideration of the relative modular classes \cite{KW}.

\section{Preliminaries: differential operators on graded commutative
  algebras}
By definition, a graded linear operator on a graded
commutative algebra $\bf A$ is of order less than or equal to $k$ if
its graded commutator with any $k+1$ left-multiplications
by elements of $\bf A$ vanishes. The graded commutator of graded
endomorphims $u$, of degree $|u|$, and $v$, of degree $|v|$,
of the graded vector space $\bf A$ is
$[u,v] = u \circ v - (-1)^{|u||v|}  v \circ u$.
Let $1$ denote the unit element of $\bf A$, and
let $\ell_a$ denote left-multiplication by $a \in \bf A$.
For $u $ a graded linear operator on $\bf A$, and $k=1$, $2$ and $3$,
we consider the 
operators $\Phi^k_u : A^{\otimes k} \to A$, defined in \cite{A}.
For $a$, $b$ and $c$ in $\bf A$, 
$$
\Phi^1_u(a)= u(a)-u(1)a,
$$
$$
\Phi^2_u(a)(b)= \Phi^1_u(ab) - \Phi^1_u(a)b - (-1)^{|a||u|} a\Phi^1_u(b),
$$
$$
\Phi^3_u(a,b)(c)= \Phi^2_u(a)(bc) - \Phi^2_u(a)(b)c - (-1)^{(|a|+|u|)|b|}
b \, \Phi^2_u(a)(c).
$$
It is easy to prove the following propositions.

\noindent
$\bullet$
$u$ is of order $0$ if and only if $\Phi^1_u = 0$.

\noindent $\bullet$ $u$ is of order $\leqslant 1$ if and only if $\Phi^2_u = 0$.

In fact,
for all $a \in {\bf A},
\Phi^2 _u(a)= \Phi^1_{[u, {\ell}_a]}$.

\noindent $\bullet$
$u$ is of order $\leqslant 1$ if and only if $\Phi^1_u$ is of order
$\leqslant 1$.

\noindent $\bullet$
A differential operator
$u$ of order $\leqslant 1$ is a derivation if and only if $u(1)=0$.

\noindent $\bullet$
$u$ is of order $\leqslant 2$ if and only if $\Phi^3_u = 0$.

\noindent $\bullet$
$u$ is of order $\leqslant 2$ if and only if $\Phi^2_u(a)$ is a derivation, for all $a \in A$.

In fact, $u$ is of order $\leqslant 2$ if and only if $[u,\ell_a]$ is of
order $\leqslant 1$ for all
$a$. This condition is equivalent to $\Phi^2_u(a)$ is of order
$\leqslant 1$ for all $a$, because $\Phi^2_u(a)$ and $[u,{\ell}_a]$ differ by
left-multiplication by $\Phi^1_u(a)$, an operator of order $0$.
Since $\Phi^2_u(a)(1)=0$, the operator
$\Phi^2_u(a)$ is of order $\leqslant 1$ if and only if
it is a derivation.

We remark that the expression $\Phi^2_u(a)(b)$ is skew-symmetric (in the
graded sense) in $a$ and $b$. More precisely,
\begin{equation}\label{skew}
(-1)^{|b|}\Phi^2_u (b)(a) = -
(-1)^{(|a|+1)(|b|+1)}(-1)^{|a|}\Phi^2_u(a)(b).
\end{equation}

\section{Bivectors and $3$-forms}\label{section3}

We shall make use of several algebraic constructions which we now describe.

\subsection{Conventions}
Hereafter $A$ is a vector bundle (later, a Lie algebroid) with base
$M$. By convention, we shall
call sections of $A$ vector fields or vectors. More generally, for $p$
and $q$ positive integers, we call
sections of $\wedge^p A$, $p$-vectors and, similarly, we call
sections of $\wedge^{q} (A^*)$, $q$-forms. The pairing of a $p$-vector
and a $p$-form is denoted $<~,~>$.

Let $i_X$ be the
interior product of forms by the vector $X$, which is a derivation of $\Gamma
(\wedge^{\bullet} A^*)$.
More generally, for vectors $X_1$ and
$X_2$, set
$$
i_{X_1 \wedge X_2}= i_{X_1} \circ i_{X_2} \ ,
$$
and define inductively the interior product of forms by a multivector.
By definition, the interior product of an $r$-form ($r$ a positive
integer) by a section of $\wedge^{q} A^* \otimes \wedge^p A$ vanishes if $p > r$ and, for $p \leqslant r$, satisfies
$$i _{\xi_1 \wedge \ldots \wedge \xi_q \otimes X_1
  \wedge \ldots \wedge X_p} (\alpha_1 \wedge \ldots \wedge\alpha_r)
 =  \xi_1 \wedge \ldots \wedge \xi_q \wedge i_{X_1
  \wedge \ldots \wedge X_p} (\alpha_1 \wedge \ldots \wedge\alpha_r) \ .
$$
Interior products by elements of $\Gamma
(\wedge^{q} A^* \otimes \wedge^p A)$
are operators of order $p$ on $\Gamma(\wedge^{\bullet} A^*)$.
The interior product of multivectors by multivector-valued forms is
similarly defined.

Given a bivector $\pi$, the vector bundle morphism $\pi^{{\sharp}}$ from
$A^*$ to $A$ is defined by
$
<\beta,\pi^{\sharp}\alpha> = \pi(\alpha, \beta)$,
for $1$-forms $\alpha$ and $\beta$. Thus
$$
i_{\pi}(\alpha\wedge\beta) =- \pi(\alpha, \beta) = <\alpha,
\pi^{\sharp}\beta > \ .
$$
Let $(e_k)$ and $(\epsilon^k)$, $1\leqslant k \leqslant N$, where $N$ is
the rank of $A$,
be dual local bases of sections of $A$ and $A^*$, respectively. Then
$$
\pi
= \frac{1}{2}  e_k \wedge
\pi^{\sharp}(\epsilon^k) =
- \frac{1}{2}
\pi^{\sharp}(\epsilon^{\ell}) \wedge e_{\ell} \ .
$$
Here and below, we use the Einstein summation convention.
The components, $\pi^{k \ell}$, of $\pi$
are defined by
$\pi= \frac{1}{2}  \pi^{k{\ell}} e_k \wedge e_{\ell}
$,
and they satisfy $\pi^{\sharp}(\epsilon^k)=  \pi^{k{\ell}}e_{\ell}$.

Let $\pi$ be a section of $\wedge ^2
A$ and $\psi$ a section of $\wedge ^3 A^*$.
We define $\psi^{(1)} $ by
$$
\psi^{(1)}(\xi)(X,Y)=
\psi(\pi^{\sharp}\xi, X,Y) \ ,
$$
for $\xi \in \Gamma(A^*)$, $X$ and $Y \in \Gamma A$.
Thus $\psi^{(1)}$ is both a vector-valued $2$-form and a
$2$-form-valued vector on $A$.
We further define $\psi^{(2)}$ by
$$
\psi^{(2)}(\xi,\eta)(X)=
\psi(\pi ^{{\sharp}}\xi, \pi ^{{\sharp}}\eta,X) \ ,
$$
for $\xi$ and $\eta$ $\in \Gamma(A^*)$, and  $X \in \Gamma A$.
Thus $\psi^{(2)}$
is both a bivector-valued $1$-form and a $1$-form-valued bivector
on $A$. In components, setting
$\psi = \frac{1}{6}
\psi_{k{\ell}m} \epsilon^k\wedge \epsilon^{\ell} \wedge \epsilon^m$, we find
\begin{equation}
\psi^{(1)}=\frac{1}{2}  \pi^{kp}\psi_{p{\ell}m} \, \epsilon^{\ell}
\wedge \epsilon^m
\otimes e_k \ ,
\end{equation}
and
\begin{equation}\label{psi2}
\psi^{(2)} = \frac{1}{2}  \pi^{kp} \pi^{{\ell}q}\psi_{pqm} \, \epsilon^m \otimes
 e_k \wedge e_{\ell} \ .
\end{equation}

\subsection{Operators on forms and multivectors}

The following operators are naturally defined on a vector bundle
$A$ equipped with a
bivector $\pi$ and a $3$-form $\psi$.

\subsubsection{The operator $\underline
  \partial_{\pi,\psi}$}

Any form-valued bivector on $A$
acts by interior product on
the sections of $\wedge ^{\bullet} A^*$.
We define the operator $\underline \partial_{\pi,\psi}$ on
sections of $\wedge ^{\bullet} A^*$ to be the interior product
by the $1$-form-valued bivector
$\psi^{(2)}$.

\begin{lemma}
a) The operator $\underline \partial_{\pi,\psi}$ is a differential
operator of order $2$ and of degree
$-1$ on the graded algebra $\Gamma(\wedge ^{\bullet} A^*)$,
which vanishes on functions and on $1$-forms.

b) For $q \geqslant 2$, and for all $\alpha_1, \dots,
\alpha_q \in \Gamma(A^*)$,
\begin{equation}\label{eq:tilde}
\underline \partial_{\pi,\psi}
(\alpha_1 \wedge \dots \wedge \alpha_q)= \sum_{1\leqslant k < \ell
  \leqslant q}
(-1)^{k+l}
\psi^{(2)}(\alpha_k ,\alpha_{\ell}) \wedge \alpha_1 \wedge \dots
\wedge \widehat{\alpha_k} \wedge \dots \wedge \widehat{\alpha_{\ell}}
\wedge \dots  \wedge \alpha_q \ .
\end{equation}
\end{lemma}
\noindent(A caret over a factor signifies that the factor is missing.)

\begin{proof}
a) We know that, for $p \in
\mathbb N$, the interior product by a form-valued $p$-vector is a
differential operator of order $p$, thus part a) follows.

b) Equation \eqref{eq:tilde} follows from the definition of the
interior product.
\end{proof}

\subsubsection{The operator $\underline d_{\pi,\psi}$}\label{dunderline}

We now consider $\psi^{(2)}$ as a bivector-valued
$1$-form on $A$, we let it act by interior product on the sections
of $\wedge ^{\bullet} A$, and we denote this operator by
$\underline{ d}_{\pi, \psi}$.

\begin{lemma}
a) The operator $\underline { d}_{\pi,\psi}$ is a derivation of
degree $+1$ of the graded algebra $\Gamma(\wedge ^{\bullet} A)$, which
vanishes on functions.

b) For $p \geqslant 1$ and for all $X_1, \dots, X_p \in \Gamma A$,
\begin{equation}\label{eq:overline}
\underline d_{\pi, \psi}(X_1 \wedge \dots \wedge X_p) =\sum_{k=1}^p
(-1)^{k+1} \psi^{(2)}(X_k) \wedge X_1 \wedge \dots
\wedge \widehat{X_k}\wedge \dots  \wedge X_p \  .
\end{equation}
\end{lemma}

\begin{proof} a) Since
$\underline{d}_{\pi, \psi}$ is a differential operator
of order $1$ and it
vanishes on functions, it is a derivation.

b) Equation \eqref{eq:overline} follows from the definition of the
interior product.
\end{proof}

\subsubsection{The operator $\delta _{\pi,\psi}$}
We define the operator $\delta_{\pi, \psi}$ on sections of
$\wedge^{\bullet}A^*$ to be
the interior product by the $2$-form-valued vector $\psi^{(1)}$.

\begin{lemma}
a) The operator $\delta_{\pi,\psi}$ is a derivation of
degree $+1$ of the graded algebra $\Gamma(\wedge ^{\bullet} A^*)$, which
vanishes on functions.

b) For $q \geqslant 1$ and for all $\alpha_1, \dots, \alpha_q \in \Gamma (A^*)$,
\begin{equation}\label{eq:hatpsi}
\delta_{\pi, \psi}(\alpha_1 \wedge \dots \wedge \alpha_q) =\sum_{k=1}^q
(-1)^{k+1}
\psi^{(1)}(\alpha_k) \wedge \alpha_1 \wedge \dots
\wedge \widehat{\alpha_k}  \wedge \dots \wedge \alpha_q \  .
\end{equation}
\end{lemma}

\begin{proof} a) Since
$\delta_{\pi , \psi}$ is the interior product by a 2-form-valued
vector, it is a derivation.

b)  Equation \eqref{eq:hatpsi} follows from the definition of the
interior product.
\end{proof}

\begin{remark}
Similarly, an operator of degree $-1$ on multivectors
can be defined
as the interior product by the vector-valued $2$-form $\psi^{(1)}$.
This operator and all three operators defined above
are $C^{\infty}(M)$-linear and can therefore
be defined pointwise.
\end{remark}

\subsection{A vector field}
Since $i_{\pi}\psi$ is a section of $A^*$, $\pi^{{\sharp}}i_{\pi}\psi$ is
a section of $A$.
We set
\begin{equation}\label{eq:defY}
    Y_{\pi,\psi}=  \pi^{{{\sharp}}} ( i_\pi \psi )  .
\end{equation}

\begin{proposition}\label{vector}
For any section $\alpha$ of $A^*$,
\begin{equation} \label{eq:y1} <\alpha,Y_{\pi,\psi}> =
 - i_{\pi^{\sharp} (\alpha) \wedge \pi} \psi =  \frac{1}{2}
   \, {\rm Tr} \Psi_{\alpha},
\end{equation}
where ${\rm Tr} \Psi_{\alpha}$ is the trace of
the endomorphism of $A^*$ defined by
\begin{equation} \label{eq:y4}
\Psi_{\alpha}(\beta)= \psi^{(2)}(\alpha, \beta),
\end{equation}
for each section $\beta$ of $A^*$.
\end{proposition}

\begin{proof}
From the definition of $Y_{\pi, \psi}$, using the skew-symmetry of
$\pi^{{\sharp}}$, we obtain
$$<\alpha, Y_{\pi, \psi}> = <\alpha, \pi^{{\sharp}}(i_{\pi}\psi)> =
- <i_{\pi}\psi, \pi^{{\sharp}}\alpha>
$$
which is indeed equal to  $- i_{\pi^{\sharp} (\alpha) \wedge \pi} \psi$,
while the trace of $\Psi_{\alpha}$ is
$$<\psi^{(2)}(\alpha, \epsilon^k),e_k > = \psi(\pi^{{\sharp}}\alpha,
\pi^{{\sharp}} \epsilon^k, e_k)
 = <\psi, \pi^{{\sharp}}\alpha \wedge \pi^{{\sharp}} \epsilon^k \wedge e_k >,
$$
where $e_k$ and $\epsilon^k$ are dual local bases of sections of $A$
and $A^*$.
The conclusion follows from the relation
$\pi = - \frac{1}{2} \, \pi^{{\sharp}} \epsilon^k \wedge e_k$.
\end{proof}

\begin{proposition}\label{commutator}
The operators on sections of $\wedge
 ^{\bullet}A^*$, $\delta_{\pi, \psi}$,
$\underline \partial _{\pi, \psi}$ and
$i_{Y_{\pi,\psi}}$ are related by
 \begin{equation}\label{eq:self}
[i_{\pi},\delta_{\pi, \psi} ]=
2 \underline \partial _{\pi, \psi}- i_{Y_{\pi,\psi}} \ .
\end{equation}
\end{proposition}

\begin{proof}
Since $i_{\pi}$ is of order 2 and of degree $-2$, and $\delta_{\pi,
  \psi}= i_{\psi^{(1)}}$
is of order 1 and of degree 1, their commutator is of order
  $\leqslant 2$
  and degree $-1$. Introducing the big bracket as in \cite{R}, we know
  that the term
  of order 2 in the commutator is the interior product by the big
  bracket $\{\pi,\psi ^{(1)}\}$ of $\pi$ and $\psi ^{(1)}$
(see \cite{yksderived}). A computation
  shows that $\psi^{(2)} = \frac{1}{2} \{\pi, \psi^{(1)}\}$, therefore
the term of order 2 is $2 i_{\psi^{(2)}}= 2 \underline \partial
_{\pi, \psi}$. If $\alpha$ is a 1-form, then $[i_{\pi},\delta_{\pi,
\psi} ]\alpha = i_{\pi} i_{\pi^{{\sharp}}\alpha}\psi$, while $$-
i_{Y_{\pi,\psi}}\alpha= - <\pi^{{\sharp}}(i_{\pi}\psi),\alpha>  =
<i_{\pi}\psi,\pi^{{\sharp}}\alpha>  = i_{\pi^{{\sharp}}\alpha}
i_{\pi} \psi.
$$
Thus equation \eqref{eq:self} is satisfied for $1$-forms.
It follows that the term of order 1 in the
commutator is $- i_{Y_{\pi,\psi}}$.
\end{proof}

\section{Lie algebroids with a twisted Poisson structure}
\label{section4}

\subsection{Lie algebroids}
We now assume that the vector bundle $A$ is a Lie algebroid over
the base manifold $M$, with anchor
$\rho$. We recall that $\rho$ is a Lie algebroid morphism from $A$ to
$TM$. We denote
the Lie bracket of
sections of $A$ and the Gerstenhaber bracket on the
graded commutative algebra, $\Gamma(\wedge^{\bullet} A)$, obtained by
extending it as a graded biderivation,
by the same symbol, $[~,~]$.
We recall that by definition $[a,.]$ is a derivation of degree $|a|-1$
of $\Gamma(\wedge^{\bullet} A)$, where $|a|$ is the degree of $a$. For
$a \in \Gamma A$ and $f \in C^{\infty}(M)$, $[a,f]= \rho(a) \cdot f$
and, for all
$a_1,\dots,a_q,b_1,\dots,b_r \in \Gamma A$, $q \geqslant 1$, $r
\geqslant 1$,
\begin{equation}\label{eq:extension}
 [a_1 \wedge \dots \wedge a_q, b_1
  \wedge \dots \wedge b_r ] \!=\! \! \sum_{k=1}^{q} \! \sum_{ {\ell}
= 1}^r
 (-1)^{k+{\ell}} [a_k,b_{\ell}]
\wedge a_1 \wedge \dots \wedge \widehat{a_{k}} \wedge  \dots  \wedge a_{q}
\wedge b_1 \wedge \dots \wedge \widehat{b_{\ell}} \wedge  \dots
\wedge b_r .
\end{equation}
We shall also
consider the differential $d_A$ on $\Gamma(\wedge^{\bullet} A^*)$
which is such that, for $f \in C^{\infty}(M)$,
$d_A f (a) = \rho(a) \cdot f$, for all $a \in \Gamma A$, and for a
$q$-form $\alpha$, $q \geqslant 1$,
$$
\displaylines{
\qquad (d_A\alpha)(a_0,...,a_{q})= \sum_{0 \leqslant k \leqslant q} (-1)^{k} \rho(a_k) \cdot
(\alpha(a_0,\dots,\widehat{a_k},\dots,a_{q})) \hfill \cr
\hfill \smash +
\sum_{0 \leqslant k < {\ell} \leqslant q} (-1)^{k+{\ell}}
\alpha([a_k,a_{\ell}],a_0,\dots,\widehat{a_k},\dots,
\widehat{a_{\ell}},\dots,a_{q}) \ , \qquad \cr}
$$
for all $a_0, \ldots, a_{q} \in \Gamma A$.
The Lie derivation of forms by a section $X$ of $A$ is the operator
${\mathcal L}^A_X = [i_X, d_A]$.
When $A$ is $TM$ with the Lie bracket of vector fields, the
differential $d_A$ is the de Rham differential of forms, and the Lie
derivation coincides with the usual notion.

\subsection{Twisted Poisson structures}

By definition, $(A, \pi, \psi)$ is a
{\it Lie algebroid with a twisted Poisson structure} if
$\pi$ is a section of $\wedge ^2
A$ and $\psi$ is a $d_A$-closed section of $\wedge ^3 A^*$ such that
\begin{equation}\label{3forme}
 \frac{1}{2}[\pi,\pi] = (\wedge^3 \pi^{{\sharp}}) \psi .
\end{equation}

To each function $f \in
C^{\infty}(M)$ is associated the section $H_f$ of $A$, called the
hamiltonian section with hamiltonian $f$, defined by
\begin{equation}\label{eq:Hamilton}
H_f= \pi^{{\sharp}}(d_Af) = - [\pi,f] .
\end{equation}

The bracket of two functions $f$ and $g$ is then defined as
$$
\{f,g\} = [H_f, g].
$$
This bracket is skew-symmetric and
satisfies the following modified Jacobi identity,
for all $f,g,h \in C^{\infty}(M)$,
\begin{equation} \label{eq:Jacobi}
\{\{f,g\},h\} + \{\{g,h\},f\} +
\{\{h,f\},g\}
= \psi(H_f, H_g, H_h)  \ .
\end{equation}
This relation is equivalent to
$$
H_{\{f,g\}} =
[H_f, H_g] + \psi^{(1)}(H_f,H_g).
$$

When a twisted Poisson structure is defined on the Lie algebroid $TM$,
the manifold $M$ is called a {\it twisted Poisson manifold}.
The following results were proved by \v Severa and Weinstein \cite{severa}
in the case of $A=TM$ and extend to the case of any Lie
algebroid $A$. See \cite{R}   for the case of Lie algebroids.
\begin{theorem}\label{dualalgd}
Let $(A, \pi, \psi)$ be a Lie algebroid with a twisted Poisson structure.
Then $A^*$ is a Lie algebroid with anchor $\rho \circ \pi^{{\sharp}}$, where
$\rho$ is the anchor of $A$, and
Lie bracket of sections, $\alpha$ and $\beta$, of $A^*$,
\begin{equation}\label{bracket}
  [\alpha,\beta]_{\pi,\psi} = [\alpha,\beta]_{\pi}  +
\psi^{(2)}(\alpha, \beta),
\end{equation}
where $[~,~]_\pi$ is defined by
\begin{equation}\label{crochetpi}
[\alpha,\beta]_{\pi} =  {\mathcal L}_{\pi^{{\sharp}} \alpha}\beta -
{\mathcal L}_{\pi^{{\sharp}} \beta}\alpha
 - d_A(\pi (\alpha, \beta)).
\end{equation}
The differential of the Lie algebroid $A^*$ is
\begin{equation}\label{diffpipsi}
d_{\pi,\psi} = d_{\pi}
+ \underline d_{\pi, \psi} ,
\end{equation}
where, for $X \in \Gamma (\wedge^{\bullet} A)$, $d_{\pi} X = [\pi,
X]$, and $\underline d_{\pi, \psi}$ is defined in \ref{dunderline}.
\end{theorem}

The map $\pi ^{{\sharp}}$ satisfies
$$
\pi ^{{\sharp}} [\alpha, \beta]_{\pi, \psi} = [ \pi ^{{\sharp}} \alpha,  \pi ^{{\sharp}}  \beta] .
$$

The case where $\psi =0$ is that of a Lie algebroid with a Poisson
structure, i.e., bivector $\pi$ such that $[\pi, \pi]=0$, and
the pair $(A,A^*)$ is also called
a {\it triangular Lie bialgebroid}.
If moreover $A=TM$, we recover the case of a
Poisson manifold.
\begin{remark}
When $\pi$ is an arbitrary bivector and $\psi$ a $3$-form,
one can still define a bracket $[~,~]_{\pi,\psi}$, which does
not in general satisfy the Jacobi identity, and a derivation,
$d_{\pi,\psi}$, which
is not in general of square zero.
\end{remark}
\begin{proposition}
Given a bivector $\pi$ and a 3-form $\psi$, the vector field
$Y_{\pi, \psi}$ defined by
equation (\ref{eq:defY}) satisfies
$$
(d_{\pi,\psi}Y_{\pi,\psi})(\alpha,\beta)
$$
$$ = - i_\pi
d_A(i_{\pi^{{\sharp}}\alpha
\wedge \pi^{{\sharp}}\beta}\psi) + < (\wedge
^3\pi^{{\sharp}}) \psi , d_A (\alpha \wedge \beta) > + < d_A \psi ,
\pi \wedge {\pi^{{\sharp}}\alpha \wedge \pi^{{\sharp}}\beta} > ,
$$
for all $\alpha$ and $\beta \in \Gamma(A^*)$.
\end{proposition}
\begin{proof}
We shall make use of the fact (see, e.g., \cite{yksm90}) that
$\wedge^{\bullet} \pi^{{\sharp}}$ is a chain map, that is to say, for any positive
integer $q$,
\begin{equation}\label{chain}
d_{\pi,\psi} \circ \wedge^{q}\pi^{\sharp} = - \wedge^{q+1}\pi^{\sharp} \circ d_A .
\end{equation}
Thus
$$ d_{\pi,\psi} Y_{\pi,\psi} = d_{\pi,\psi}\pi^{\sharp} (i_{\pi} \psi )
= - (\wedge^2 \pi^{\sharp})  d_A i_\pi  \psi, $$
which implies
$$(d_{\pi,\psi} Y_{\pi,\psi}) (\alpha, \beta) =
i_{\pi^{{\sharp}}\alpha \wedge \pi^{{\sharp}}\beta} d_A i_{\pi} \psi
$$
$$
= i_{\pi^{{\sharp}}\alpha \wedge \pi^{{\sharp}}\beta} [d_A, i_{\pi}] \psi
- <i_{\pi}d_A\psi, {\pi^{{\sharp}}\alpha \wedge \pi^{{\sharp}}\beta} > .
$$
The first term is
$$
 i_{\pi^{{\sharp}}\alpha \wedge \pi^{{\sharp}}\beta} [d_A, i_{\pi}] \psi =
[ i_{\pi^{{\sharp}}\alpha \wedge \pi^{{\sharp}}\beta}, [d_A, i_{\pi}]] \psi
+[d_A, i_{\pi}]
 i_{\pi^{{\sharp}}\alpha \wedge \pi^{{\sharp}}\beta}  \psi .
$$
The analogue of the Cartan relation for Lie algebroids (see, e.g.,
\cite{yksderived})
shows that
$$[ i_{\pi^{{\sharp}}\alpha \wedge \pi^{{\sharp}}\beta}, [d_A, i_{\pi}]]
= i_{[\pi^{{\sharp}}\alpha \wedge \pi^{{\sharp}}\beta, \pi]} ,
$$
where the bracket on the left-hand side is the graded commutator,
while the bracket on the right-hand side is the Gerstenhaber bracket
of sections of $\wedge^{\bullet}A$.
To conclude the proof, we express $  < [\pi^{\sharp} \alpha   \wedge  \pi^{{\sharp}}
\beta, \pi],\psi > $ in terms of $d_A(\alpha \wedge \beta)$.
By equation \eqref{chain},
$$  [\pi^{{\sharp}} \alpha   \wedge  \pi^{{\sharp}} \beta, \pi]
   = d_{\pi}  (  \pi^{\sharp} \alpha   \wedge  \pi^{\sharp} \beta  ) $$
$$ =d_{\pi,\psi}(\pi^{\sharp} \alpha \wedge
\pi^{\sharp} \beta) - {\underline d}_{\pi,\psi}
(\pi^{\sharp} \alpha \wedge  \pi^{\sharp}
\beta)=    -   ( \wedge^3 \pi^{\sharp} )   d_A (  \alpha \wedge \beta
) - {\underline d}_{\pi,\psi}
(\pi^{\sharp} \alpha \wedge  \pi^{\sharp}  \beta)  .$$
By equation \eqref{eq:overline},
$$
<{\underline d}_{\pi,\psi} (\pi^{\sharp} \alpha \wedge \pi^{\sharp} \beta) ,\psi >
$$
$$
=  <   (\wedge^2 \pi^{{\sharp}} )(i_{\pi^{{\sharp}} \alpha }\psi) \wedge
  \pi^{{\sharp}} \beta , \psi > - <(\wedge^2 \pi^{{\sharp}}) (i_{\pi^{{\sharp}}
    \beta }\psi) \wedge
\pi^{{\sharp}} \alpha , \psi > =0 ,
$$
since the operator $\wedge^2 \pi^{\sharp}$
is symmetric.
Therefore,
$$  < [\pi^{\sharp} \alpha   \wedge  \pi^{\sharp} \beta, \pi],\psi >  =  -    <    ( \wedge^3 \pi^{\sharp} )   d_A (  \alpha \wedge \beta
),\psi >. $$
The proposition follows.
\end{proof}
\begin{corollary}\label{corollary}
If $(A, \pi,\psi)$ is a Lie algebroid with a twisted
Poisson structure, the $d_{\pi,\psi}$-coboundary of $Y_{\pi,
    \psi}$
satisfies, for all $\alpha$ and $\beta \in \Gamma(A^*)$,
\begin{equation} \label{coboundary}
 (d_{\pi,\psi}Y_{\pi,\psi})(\alpha,\beta)  = - i_\pi
d_A(i_{\pi^{{\sharp}}\alpha \wedge \pi^{{\sharp}}\beta}\psi) + \frac{1}{2}
<[\pi,\pi], d_A (\alpha \wedge \beta) > .
\end{equation}
\end{corollary}
We shall make use of this formula in the proof of Theorem
\ref{generatorsquare0}.

\section{Generators of the Gerstenhaber algebra of a twisted Poisson
structure}
\label{section5}
By definition, a {\it generator} of a Gerstenhaber algebra,  $(\bf A,
[~,~]_{\bf A})$, is an operator $u$ on $\bf A$ such that
\begin{equation}\label{generator}
[a,b]_{\bf A} = (-1)^{|a|}(u(ab)-u(a)b - (-1)^{|a|}a u(b)),
\end{equation}
for all $a$ and $b \in {\bf A}$.

It is clear that any two generators of a Gerstenhaber algebra differ by a
derivation of the underlying graded commutative algebra.

\begin{theorem}\label{agenerator}
The operator $\partial_{\pi}+\underline \partial_{\pi,\psi}$, where
$\partial_{\pi} = [d_A,i_{\pi} ]$,
is a generator of the
Gerstenhaber algebra $(\Gamma(\wedge ^{\bullet} A^*), [~,~]_{\pi,\psi})$
associated to
the Lie algebroid $A$ with the twisted Poisson structure $(\pi, \psi)$.
\end{theorem}

\begin{proof} By definition of the Lie algebroid bracket of $A^*$, for
$1$-forms, $\alpha$ and $\beta$,
$$
[\alpha,\beta]_{\pi,\psi} =[\alpha,\beta]_{\pi} + \psi^{(2)}(\alpha,
\beta).
$$
The differential operator, $\partial_{\pi}$, is of
order 2 and of degree $-1$. Therefore, for each $\alpha \in
\Gamma(\wedge^{\bullet} A^*)$, the map
$$\beta \mapsto (-1)^{|\alpha|}\Phi^2_{\partial_{\pi}}(\alpha)(\beta)
= (-1)^{|\alpha|} (\partial_{\pi}(\alpha \wedge \beta) -
\partial_{\pi} \alpha \wedge \beta -
(-1)^{|\alpha|} \alpha \wedge  \partial_{\pi}\beta)
$$
is a derivation of degree $|\alpha|-1$. Since it satisfies the
condition of skew-symmetry \eqref{skew}, it is enough to show that it
coincides with $[\alpha,\beta]_\pi$ when $\alpha$ is of degree 1 and
$\beta$ is of degree 0 or 1.
If $\beta$ is a 0-form, $f$, then
$$(-1)^{|\alpha|}\Phi^2_{\partial_{\pi}}(\alpha)(\beta) = 
(-1)^{|\alpha|}(d_Af\wedge i_{\pi} \alpha -i_{\pi}(d_Af\wedge \alpha))
\ .
$$
When $\alpha$ is a 1-form, this expression is equal to
$i_{\pi}(d_Af\wedge\alpha)= <d_Af, \pi^{{\sharp}}\alpha>$, which is
$[\alpha,f]_{\pi}$ by definition.
When $\alpha$ and $\beta$ are 1-forms,
$$(-1)^{|\alpha|}
\Phi^2_{\partial_{\pi}}(\alpha)(\beta)
$$
$$
= d_A(\pi(\alpha,\beta))+i_{\pi}(d_A\alpha\wedge\beta)
-i_{\pi}(\alpha\wedge d_A\beta)-(i_{\pi}d_A\alpha)\beta
+\alpha (i_{\pi}d_A\beta) \ ,
$$
while
$$
[\alpha,\beta]_\pi = i_{\pi^{{\sharp}}\alpha}d_A\beta
-i_{\pi^{{\sharp}}\beta}d_A\alpha +d_A(\pi(\alpha,\beta)) \ ,
$$
and both expressions coincide since
$i_{\pi}(d_A\alpha \wedge \beta) -(i_{\pi} d_A\alpha) \beta = 
- i_{\pi^{{\sharp}}\beta}d_A\alpha$,
as can be easily shown.

We now consider the bilinear map $(\alpha, \beta) \mapsto (-1)^{|\alpha|}
\Phi^2_{\underline \partial_{\pi,\psi}}(\alpha)(\beta)$ which coincides
with $\psi^{(2)}$ for $\alpha$ and $\beta$ of degree
$1$ and vanishes if $\alpha$ or $\beta$ is of degree 0. By
\eqref{skew}, it is skew-symmetric in the graded sense and, for each
$\alpha$, it defines a derivation of degree $|\alpha|-1$, since
${\underline \partial}_{\pi,\psi}$ is a differential operator of order 2.
Therefore this bilinear map is the required extension of $\psi^{(2)}$.
\end{proof}

But the square of the generator
$\partial_{\pi}+\underline \partial_{\pi,\psi}$ does not vanish in
general.

\begin{remark}
When $\pi$ is an arbitrary bivector and $\psi$ a $3$-form,
the extension of
bracket $[~,~]_{\pi,\psi}$ as a
biderivation on $\Gamma(\wedge^{\bullet}A^*)$ is not in general
a Gerstenhaber algebra bracket. However the proof of the preceding theorem
shows that
the operator $\partial_{\pi}+\underline \partial_{\pi,\psi}$
is a ``generator'' of this bracket, in the sense
that it satisfies relation \eqref{generator}.
\end{remark}

\begin{lemma}\label{lemmacocycle}
Let $\partial$ be a generator of the Gerstenhaber algebra $\Gamma
(\wedge^{\bullet}E)$ of a Lie algebroid $E$, and
let $U$ be a section of $E^*$. Then
$\partial +i_U$ is a generator of square zero of
$\Gamma  (\wedge^{\bullet}E)$ if and only if
$\partial ^2= i_{d_E U}$.
\end{lemma}
\begin{proof}
The generalization of formula (2.4)
in \cite{Ko} (see, e.g., \cite{yks2000}) implies
$$(\partial +i_U)^2 = \partial ^2 + [\partial, i_U] = \partial^2 -
i_{d_E U},
$$
which proves the claim.
\end{proof}

\begin{theorem}\label{generatorsquare0}
The operator $\partial_{\pi}+\underline
\partial_{\pi,\psi}+i_{Y_{\pi,\psi}}$
is a generator of square zero of the
Gerstenhaber algebra $(\Gamma(\wedge ^{\bullet} A^*), [~,~]_{\pi,\psi})$
associated to
the Lie algebroid $A$ with the twisted Poisson structure $(\pi, \psi)$.
\end{theorem}

\begin{proof}
For 1-forms $\alpha$ and $\beta$, compute
$$(\partial_{\pi}+\underline \partial_{\pi,\psi})(\alpha \wedge \beta) =
- d_A (\pi(\alpha, \beta)) - i_{\pi} d_A (\alpha \wedge \beta) +
{i}_{\pi^{{\sharp}}\alpha \wedge \pi^{{\sharp}}\beta} \psi  \ .$$
Since $ \underline \partial_{\pi,\psi}  $ vanishes on $1$-forms,
and since $ \partial_{\pi} d_A (\pi(\alpha, \beta))= 0$,
$$
(\partial_{\pi}+\underline \partial_{\pi,\psi})^2 (\alpha \wedge \beta)
 = i_{\pi} d_A i_{\pi} (d_A(\alpha \wedge \beta))
- i_{\pi} d_A ({i}_{\pi^{{\sharp}}\alpha
\wedge \pi^{{\sharp}}\beta} \psi)  \ .$$
For any closed 3-form $\tau$, $[[i_{\pi}, d_A], i_{\pi}] \tau = 2
i_{\pi} d_A i_{\pi} \tau$.
Since $[[i_{\pi},d_A],i_{\pi}]=i_{[\pi,\pi]}$, where the
bracket on the right-hand side is the
Gerstenhaber bracket of multivectors, we obtain
$$i_{\pi} d_A i_{\pi} \tau = \frac{1}{2} [[i_{\pi},d_A],i_{\pi}] \tau =
\frac{1}{2}i_{[\pi,\pi]} \tau .
$$
Therefore,
in view of Corollary \ref{corollary}, when $\pi$ and $\psi$ satisfy
the equations of a twisted Poisson structure,
$$
(\partial_{\pi}+\underline \partial_{\pi,\psi})^2 (\alpha \wedge \beta)
=i_{d_{\pi,\psi}Y_{\pi,\psi}}(\alpha \wedge \beta) .
$$
The theorem then follows from Theorem \ref{agenerator}
and the preceding lemma applied to the Lie
algebroid $E = A^*$.
\end{proof}

We recall that a Gerstenhaber algebra equipped with a generator of
square zero is called a {\it Batalin-Vilkovisky algebra} or, for
short, a {\it BV-algebra}.
We can reformulate the preceding theorem as follows.

If $(A,\pi,\psi)$ is a Lie algebroid with a twisted Poisson structure,
{\it the algebra $(\Gamma(\wedge^{\bullet}A^*), [~,~]_{\pi,\psi})$
is a BV-algebra, with generator
$\partial_{\pi} +
\underline \partial _{\pi, \psi}+ i_{Y_{\pi,\psi}}$}.

\section{Generators and the modular class}\label{sectionZ}

For simplicity, we shall assume that the vector bundle $A$ is
orientable, i.e., admits a
nowhere-vanishing form of top degree,
$\lambda \in \Gamma(\wedge^{N}A^*)$, where
$N$ is the rank of $A$, which we call a volume form.
If the bundle is non-orientable, densities
must be used instead of volume forms, and the proofs need not be changed.

\subsection{A generator of square zero}

The $N$-form $\lambda$ defines an isomorphism $*_{\lambda}$
of vector bundles
from $\wedge^{\bullet} A$ to $\wedge^{\bullet} A^*$ by
$*_{\lambda}V = i_V\lambda$, which induces
a map on sections denoted in the same way. If $V$ is a
$p$-vector, then $*_{\lambda}V$ is an $(N-p)$-form.

Given a twisted Poisson structure on the Lie algebroid $A$,
consider the operator
$$
\partial_{\pi, \psi,\lambda} =-*_{\lambda} d_{\pi,\psi} \, *^{-1}_{\lambda},
$$
where $d_{\pi,\psi}$ is the differential \eqref{diffpipsi}.

\begin{proposition}\label{square0}
The operator $\partial_{\pi, \psi,\lambda}$
is a generator of square zero of the
Gerstenhaber algebra $(\Gamma(\wedge ^{\bullet} A^*), [~,~]_{\pi,\psi})$
associated to
the Lie algebroid $A$ with the twisted Poisson structure $(\pi, \psi)$.
\end{proposition}

\begin{proof}
The square
of $\partial_{\pi, \psi, \lambda}$ vanishes since $d_{\pi,\psi}$ is
a differential. It is well-known that
conjugating the differential of the Lie algebroid by $*_{\lambda}$
yields the opposite of a
generator of the Gerstenhaber bracket. See, e.g., \cite{yks2000}.
\end{proof}

\subsection{Properties of sections $X_{\pi, \lambda}$ and $Y_{\pi, \psi}$}
In order to generalize the modular vector fields of Poisson manifolds,
we shall first consider the section $X_{\pi, \lambda}$ of $A$ such
that,
for all $ \alpha \in \Gamma(A^*)$,
\begin{equation}\label{defX}
  < \alpha, X_{\pi,\lambda}> \lambda  =
{\mathcal L}^A_{ \pi^{{\sharp}} \alpha  } \lambda - (i_{\pi} d_A \alpha )
  \lambda \ ,
\end{equation}
where ${\mathcal L}^A$ is the Lie derivation.
Since the right-hand side of the previous
expression is $C^{\infty}(M)$-linear in $\alpha$, the section
$X_{\pi,\lambda}$ is well-defined.

In particular, ${\mathcal L}^A_{H_f} \lambda = <d_A f, X_{\pi,\lambda}>
\lambda$, where, as above, $H_f$ is
the hamiltonian section with hamiltonian $f \in C^{\infty}(M)$.
Thus, when $A=TM$, the vector field $X_{\pi,\lambda}$ satisfies, for
each $f \in C^{\infty}(M)$,
\begin{equation}\label{propX}
{\mathcal L}_{H_f} \lambda = (X_{\pi,\lambda} \cdot f)\lambda \ ,
\end{equation}
i.e.,
the function $X_{\pi, \lambda} \cdot f$ is the divergence of
$H_f$ with respect to the volume form~$\lambda$.

\begin{lemma}\label{lemmaX}
The section $X_{\pi,\lambda}$ satisfies
\begin{equation}\label{*X}
\partial_{\pi} \lambda = - i_{X_{\pi,\lambda}}\lambda .
\end{equation}
\end{lemma}

\begin{proof}
We shall make use of the fact that
for any $\eta \in \Gamma ( \wedge
^{\bullet} A^*)$,
$\eta \wedge i_{\pi} \lambda = i_{\pi}\eta \wedge \lambda$.
By the definition of $X_{\pi,\lambda}$,
$$
\alpha \wedge i_{X_{\pi,\lambda}} \lambda = d_A (i_{\pi^{{\sharp}}\alpha} \lambda)
- (i_{\pi} d_A \alpha) \lambda ,
$$
for each 1-form $\alpha$.
Since $\lambda$ is a form of top degree, $
i_{\pi^{{\sharp}}\alpha}\lambda =  \alpha \wedge i_{\pi} \lambda .
$
Therefore
$$
\alpha \wedge i_{X_{\pi,\lambda}} \lambda =
d_A \alpha \wedge i_{\pi} \lambda - \alpha \wedge d_A i_{\pi} \lambda
- (i_{\pi} d_A \alpha) \lambda
= - \alpha \wedge d_A i_{\pi} \lambda \ .
$$
Thus $i_{X_{\pi,\lambda}}\lambda = - d_A i_{\pi}\lambda = -\partial_{\pi}\lambda$.
\end{proof}

We now consider the section $Y_{\pi,\psi}$ defined by \eqref{eq:defY}.
\begin{lemma}\label{lemmaY}
The section $Y_{\pi, \psi}$ satisfies
\begin{equation}\label{propY}
{\underline {\partial}}_{\pi, \psi}  \lambda =
- 2 \, i_{Y_{\pi, \psi}} \lambda \ .
\end{equation}
\end{lemma}
\begin{proof}
Let $\xi \wedge Q$ be a decomposable form-valued bivector.
Then
$$i_{\xi \wedge Q} \lambda =
\xi \wedge i_Q \lambda = \varepsilon_{\xi} *_{\lambda} Q
= *_{\lambda} \, i_{\xi} Q = i_{i_{\xi}Q}\lambda \ ,
$$
where $\varepsilon_\xi$ is the left exterior product by $\xi$.
Let $(e_k)$ and $(\epsilon^k)$, $1\leqslant k \leqslant N$,
be dual local bases of sections of $A$ and $A^*$, such that
$\lambda = \epsilon^1 \wedge \dots \wedge \epsilon^N$.
Let $\sigma =  \frac{1}{2} \sigma^{k\ell}_m \epsilon^m \otimes e_k
\wedge e_{\ell}$. Then $
i_{\sigma} \lambda =
i_S\lambda$,
where $S = \sigma^{k\ell}_k e_{\ell}$ is twice the trace of
$\sigma$ with respect to the first index.
If $\sigma = \psi^{(2)}$, by formulas \eqref{psi2} and \eqref{eq:y1},
$S= - 2 Y_{\pi, \psi}$.
\end{proof}

We remark that $\pi$ and
$\psi$
need not satisfy the axioms of a twisted Poisson structure in order
for the results of these two lemmas to be valid. Lemma \ref{lemmaX}
expresses the fact that the isomorphism $*_{\lambda}$ identifies the
vector field $X_{\pi,\lambda}$ with the $(N-1)$-form $-
\partial_{\pi}\lambda= - d_A i_{\pi} \lambda$,
a property which is valid for the modular vector fields of Poisson
manifolds  \cite{yks2000} \cite{W}.

\subsection{The modular class}

Set
\begin{equation}
Z_{\pi, \psi, \lambda} = X_{\pi, \lambda} +Y_{\pi, \psi} \ ,
\end{equation}
where $Y_{\pi,\psi}$ is defined by \eqref{eq:defY}.

\begin{proposition}
The section $Z_{\pi, \psi, \lambda}$ satisfies the relation
\begin{equation}\label{difference}
\partial_{\pi, \psi,\lambda}- (\partial_{\pi} + \underline
\partial_{\pi, \psi} + i_{Y_{\pi, \psi}}) = i_{Z_{\pi, \psi, \lambda}}.
\end{equation}
\end{proposition}
\begin{proof}
Because both $\partial_{\pi, \psi,\lambda}$ and $\partial_{\pi} + \underline
\partial_{\pi, \psi} + i_{Y_{\pi, \psi}}$ are generators of the same
Gerstenhaber algebra, they differ by the interior product by a section of
$A$. It is enough to evaluate their difference on the form of top degree
$\lambda$, and the result follows from the fact that
$\partial_{\pi, \psi,\lambda} \, \lambda = 0$, together with
Lemmas \ref{lemmaX} and \ref{lemmaY}.
\end{proof}

By the general properties of \cite{yks2000},
we obtain

\begin{theorem}\label{theorem6.1}
The section
$Z_{\pi, \psi, \lambda} = X_{\pi, \lambda} +Y_{\pi, \psi}$
of $A$  is a $d_{\pi, \psi}$-cocycle. The cohomology
class of $Z_{\pi, \psi, \lambda}$ is independent of the choice of $\lambda$.
\end{theorem}

\begin{proof}
In fact, it follows from Lemma \ref{lemmacocycle} that
the difference of two generating operators whose squares
vanish is the interior product by a $1$-cocycle. So the fact that
$d_{\pi,\psi}Z_{\pi,\psi,\lambda} = 0$ is a
consequence of Theorem \ref{generatorsquare0}
and of Proposition \ref{square0}. Replacing the form of top degree $\lambda$ by
the form $f \lambda$, where $f$ is a nowhere-vanishing smooth function on the
base manifold, adds the coboundary $d_{\pi} (\ln |f|)
= d_{\pi, \psi} (\ln |f|)$ to $X_{\pi, \lambda}$, therefore
$Z_{\pi,\psi,f \lambda}$ and $Z_{\pi,\psi,\lambda}$ are cohomologous cocycles.
\end{proof}

\begin{remark}
One can prove directly that $d_{\pi,\psi}Z_{\pi,\psi,\lambda}=0$, but
we have not found a proof simpler than the one given here.
\end{remark}
\begin{definition} The section $Z_{\pi, \psi, \lambda}$ is called a
{\em modular section} (or modular vector field) of $(A,\pi,\psi)$.
The $d_{\pi, \psi}$-cohomology class of
$Z_{\pi, \psi, \lambda}$
is called the {\em modular class} of
the Lie algebroid $A$ with the twisted Poisson structure $(\pi,
\psi)$, and $A$ is called {\em unimodular} if its modular class vanishes.
When $A=TM$, the modular class of $(TM, \pi,\psi)$ is called the
{\em modular class of the twisted Poisson manifold} $(M, \pi, \psi)$.
\end{definition}

The modular classes of twisted Poisson manifolds generalize
the modular classes of Poisson manifolds.
If $\psi =0$, then $Z_{\pi,\psi, \lambda}= X_{\pi,\lambda}$, the
generator $\partial_{\pi,\psi,\lambda}$ reduces to $\partial_{\pi,
\lambda}= - *_{\lambda} d_{\pi} *^{-1}_{\lambda}$, and
relation \eqref{difference} reduces to
\begin{equation}\label{poisson}
\partial_{\pi,\lambda}- \partial_{\pi} = i_{X_{\pi, \lambda} \ ,}
\end{equation}
a relation valid for the modular vector fields
of a triangular Lie bialgebroid, in particular of a Poisson manifold
\cite{yks2000}.

\subsection{Properties of the modular sections}
We now list properties of the modular sections which
generalize the properties of the modular vector fields of Poisson
manifolds.
The modular vector fields of Poisson manifolds
and triangular Lie bialgebroids satisfy
$$
 *_{\lambda} X_{\pi,\lambda} = - \partial_{\pi}
\lambda \ ,
$$
whereas in the twisted case,
equation
\eqref{difference} implies
\begin{equation}\label{*Z}
*_{\lambda} Z_{\pi,\psi,\lambda} = - (\partial_{\pi} +
\underline{\partial}_{\pi, \psi} +i_{Y_{\pi,\psi}})\lambda \ .
\end{equation}
The modular vector fields of Poisson manifolds
and triangular Lie bialgebroids satisfy
$ *_{\lambda} X_{\pi,\lambda} = -  d_{A} *_\lambda \pi$, whereas
in the twisted case, in view of Proposition \ref{commutator},
\begin{equation}\label{*Z2}
*_{\lambda}
Z_{\pi,\psi, \lambda} =  - (d_{A}+\delta_{\pi, \psi})
*_\lambda \pi - 3 \, \underline{\partial}_{\pi,
  \psi}\lambda \ .
\end{equation}

In the untwisted case, the modular vector fields
satisfy relation \eqref{defX} for any
1-form $\alpha$,
and therefore
$$  < \alpha, X_{\pi,\lambda}> \lambda  =
{\mathcal L}^A_{ \pi^{{\sharp}} \alpha  } \lambda \,
$$
for any
$d_A$-closed 1-form. In the twisted case, adding
relations \eqref{defX} and \eqref{eq:y1}, we obtain
$$ < \alpha,Z_{\pi,\psi,\lambda}> \lambda = {\mathcal L}^A_{ \pi^{\sharp} \alpha } \lambda
- i_{\pi} ((d_A + \delta_{\pi, \psi})\alpha) \lambda \ .
$$
We have used the operator
$\delta_{\pi, \psi}$ introduced in
\eqref{eq:hatpsi}
to write the term
$i_{\pi^{\sharp}\alpha} \psi$ as $\delta_{\pi, \psi}\alpha$.
The differential $d_A$ of the Lie algebroid $A$
is twisted into the derivation
\begin{equation}\label{quasidiff}
d_{A,\pi,\psi} = d_A + \delta_{\pi, \psi}
\end{equation}
of $\Gamma (\wedge
^{\bullet}A^*)$, which is no longer of square zero, and
the modular section $Z_{\pi, \psi, \lambda}$
satisfies
\begin{equation}
 < \alpha,Z_{\pi,\psi,\lambda}> \lambda = {\mathcal L}^A_{ \pi^{\sharp} \alpha } \lambda \ ,
\end{equation}
for any $d_{A,\pi,\psi}$-closed $1$-form $\alpha$.

\subsection{The unimodular case}

When $(A,\pi,\psi)$ is
unimodular, i.e., the class of the
  modular section $Z_{\pi,
  \psi, \lambda}$ vanishes, the homology and cohomology are
isomorphic, the isomorphism being, in fact, defined at the chain level.
By definition, 
the {\em homology} $H_\bullet^{\pi,\psi} (A) $ of a Lie algebroid
with a twisted
Poisson structure, $(A,\pi,\psi)$, is
the homology of the complex
  $(\Gamma(\wedge^\bullet A^*),\partial_{\pi} + \underline
  \partial_{\pi, \psi} + i_{Y_{\pi, \psi}}  ) $.
The untwisted case, generalizing the Poisson homology of Koszul and
Brylinski, was studied by Huebschmann in the framework of
Lie-Rinehart algebras \cite{H} \cite{H2} and by Xu \cite{Xu}.
The {\em cohomology} $ H^{\bullet}_{\pi,\psi} (A^*) $
of a Lie algebroid with a twisted Poisson structure
is  the cohomology of the Lie algebroid $A^*$
defined in Theorem \ref{dualalgd},
i.e., the cohomology of the complex $ (\Gamma(\wedge^\bullet A) , d_{\pi,\psi})  $.

\begin{proposition}
If $(A,\pi,\psi)$ is unimodular, for all $k \in {\mathbb N}$,
$H_k^{\pi,\psi} (A)  \simeq H^{N-k}_{\pi,\psi} (A^*) $, where $N$ is
the rank of $A$.
\end{proposition}

This proposition follows from the results of \cite{H} and
\cite{H2}, or \cite{Xu}.
However, to make this paper self-contained, we present a proof.

\begin{proof}
Let $ \lambda$ be a
 nowhere-vanishing form of top degree, and
compare $\partial_{\pi,\psi,f\lambda} $
and $\partial_{\pi,\psi,\lambda}$, where $f$ is a nowhere-vanishing
 function on $M$.
Then, by the graded Leibniz identity, for
any multivector $X$ and any function $g $,
$$
d_{\pi}  (g X)
= [\pi, g] \wedge X +  g \ d_{\pi} X
=   -H_g \wedge X +
g \ d_{\pi} X.$$
Since $ \underline d_{\pi,\psi}$
is $C^{\infty}(M)$-linear, this relation implies
$$
d_{\pi,\psi}  (g X)
=   -H_g \wedge X +  g \ d_{\pi,\psi} X .
$$
For any $p$-form $\alpha$,
$$
 \partial_{\pi,\psi,f \lambda}(\alpha) =  -  *_{f\lambda}
 d_{\pi,\psi}    *_{f \lambda}^{-1} \alpha
=   -  *_{\lambda}  d_{\pi,\psi}  *_{\lambda}^{-1} \alpha   +  f *_{\lambda}
   ( H_{\frac{1}{f} } \wedge  *_{\lambda}^{-1} \alpha)
$$
$$ =     \partial_{\pi,\psi, \lambda} \alpha   +
f i_{ H_{\frac{1}{f}}} \alpha
= \partial_{\pi,\psi, \lambda}\alpha +
i_{H_{\ln |f|}} \alpha .
$$
We have used the fact that $i_X \circ *_{\lambda} = *_{\lambda} \circ
\varepsilon _X$,
where $\varepsilon_X$ is the left exterior product by $X$, which implies that
$*_{\lambda}
( X \wedge *_\lambda^{-1} \alpha)= i_X \alpha. $

If the modular class vanishes, there exists
$g \in C^{\infty} (M)$ such that $ Z_{\pi,\psi,\lambda}= H_g$.
Set $f = e^{-g}$. Then
  $$  \partial_{\pi} + \underline \partial_{\pi, \psi} + i_{Y_{\pi,
      \psi}} =
\partial_{\pi,\psi, \lambda} -i_{Z_{\pi,\psi,\lambda}} = 
\partial_{\pi,\psi,\lambda}  - i_{H_g} =  \partial_{\pi,\psi,f
 \lambda} \ .$$
In other words, the map $V \mapsto  *_{f \lambda} V$ is
a chain map from
  $ (\Gamma(\wedge^\bullet A) , d_{\pi,\psi})  $ to
$(\Gamma(\wedge^\bullet A^*),\partial_{\pi} + \underline \partial_{\pi
   , \psi}+ i_{Y_{\pi, \psi}}  ) $.
\end{proof}

\section{Comparison with the ELW-modular class}
\label{sectionELW}

As stated in Theorem \ref{dualalgd},
when $(A,\pi,\psi)$ is a Lie algebroid with a twisted Poisson
structure, the dual vector bundle
$A^*$ is a Lie algebroid with anchor $\rho \circ
\pi^{{\sharp}}$ and bracket
$[~,~]_{\pi,\psi}$.
Since in \cite{H2} and  \cite{ELW}, general notions for
Lie algebroids were defined, a comparison is in order.
In \cite{ELW}, Evens, Lu and Weinstein defined the characteristic class
of a Lie algebroid with a representation in a line bundle and the
modular class of a Lie
algebroid $E$ -- which we shall call the ELW-modular class of
$E$.
We shall compare what we have called the
modular class of $(A,\pi, \psi)$ to the ELW-modular class
of the Lie
algebroid
$(A^*, \rho \circ \pi^{{\sharp}}, [~,~]_{\pi,\psi})$, and conclude that in
the case of $A=TM$, the first is one half the second, a
result that is not valid in the case of a Lie algebroid in general.

\subsection{The characteristic class}
We recall the construction of \cite{ELW}.
Let $E$ be a Lie algebroid over a manifold $M$, with anchor $\rho$ and
Lie bracket of sections $[~,~]_E$, and let $D$ be a representation of $E$
on a line bundle $L$ over $M$,
$x\in \Gamma E \mapsto D_x \in {\rm {End}}(\Gamma L)$. By definition, the map
$D$ is $C^{\infty}(M)$-linear and  $D_x (f \mu)= f D_x \mu
 + ( \rho(x) \cdot f ) \, \mu$, for all $x \in \Gamma E$,
$\mu \in \Gamma L$ and $f \in C^{\infty}(M) $, and
$D_{[x,y]_E}=[D_x,D_y]$, for all $x$ and $y$ in $\Gamma E$.
The {\it characteristic class of $E$ associated to the representation $D$
on $L$} is the class of the $d_E$-cocycle $\theta_{s} \in \Gamma
(E^*)$
defined by
$$D_x s = < \theta_s , x > s,
$$
where $s$ is a nowhere-vanishing section
of $L$. If $L$ is not trivial, the class of $L$ is defined as one
half that of its square.

Assume that $\partial$ is a generating operator of the Gerstenhaber
bracket, $[~,~]_E$, of $\Gamma (\wedge^{\bullet}E)$. Set
\begin{equation}\label{rep}
D^{\partial}_x (\mu) = [x,\mu]_E - (\partial x) \mu = - x \wedge \partial \mu \ ,
\end{equation}
for $x \in \Gamma E$ and $\mu \in \Gamma (\wedge^{N}E)$. Then
 $D^{\partial}$
is a
representation of $E$ on $\wedge^{N}E$, and the
associated characteristic class is the class
of the cocycle $\xi \in \Gamma (E^*)$ such that
$$
<\xi,x> \mu = - x \wedge \partial \mu \ .
$$

If, in particular, $(A,\pi, \psi)$ is a Lie algebroid with a
twisted Poisson structure,  we can consider the Lie algebroid $E = A^*$,
with Lie bracket of sections
$[~,~]_{\pi, \psi}$, and generator $\partial =
\partial _{\pi} +
{\underline {\partial}}_{\pi,\psi} +i_{Y_{\pi, \psi}}$.
If $\lambda$ is a nowhere-vanishing section of
$\Gamma (\wedge^{N}A^*)$ (or a density in the non-orientable case),
then, by equation \eqref{*Z},
$\partial \lambda = - i_{Z_{\pi, \psi, \lambda}}\lambda$. The associated
characteristic class is the class of the $d_{\pi, \psi}$-cocycle
$\theta \in \Gamma A$ such that
$$
<\alpha, \theta > \lambda = \alpha \wedge i_{Z_{\pi, \psi, \lambda}} \lambda = <\alpha ,
Z_{\pi,
  \psi, \lambda} >  \lambda \ ,
$$
for all $\alpha \in \Gamma (A^*)$.
Therefore,
\begin{proposition}
The characteristic class of $A^*$ associated to the
representation \eqref{rep}, $\alpha \mapsto D^{\partial}_{\alpha}$,
of $A^*$ on $\wedge^N (A^*)$, where
 $\partial =
\partial _{\pi} +
{\underline {\partial}}_{\pi,\psi} +i_{Y_{\pi, \psi}}$,
coincides with the modular class of
the Lie algebroid $A$ with the twisted Poisson structure
$(\pi, \psi)$.
\end{proposition}

\subsection{The ELW-modular class of $A^*$}

In \cite{ELW}, the modular class of a Lie algebroid $E$
is defined as follows.
Let $ L^E= \wedge^{N} E \otimes \wedge^{n} T^*M $,
where $n$ is the dimension of $M$.
Define a representation $D^E$ of $E$ on $L^E$ by
\begin{equation} \label{representation}
D^E_x(\omega \otimes \mu) = [x, \omega]_E \otimes \mu + \omega
\otimes {\mathcal L}_{\rho(x)}\mu ,
\end{equation}
for $x \in \Gamma E$ and $\omega \otimes \mu \in \Gamma (L^E)$. Here $\mathcal L$ is the
Lie derivation of forms on $M$ by vector fields.

\begin{definition} The {\em ELW-modular class} of the Lie algebroid $E$ is the
characteristic class of $E$ associated to the representation $D^E$
of $E$ on $L^E$.
\end{definition}

\begin{proposition} \label{prop:evens}
The modular class
of a twisted Poisson manifold is equal to one half the ELW-modular
class of its cotangent bundle Lie algebroid.
\end{proposition}
\begin{proof}
Let $\lambda $ be a volume form on $M$, i.e., a nowhere-vanishing
section of $\wedge^n T^*M$.
By definition, the ELW-modular class of the Lie algebroid $T^*M$ is
the class of the vector field $U$ such that
\begin{equation} \label{U}
<\alpha, U > \lambda \otimes \lambda = [\alpha, \lambda]_{\pi, \psi}
\otimes \lambda + \lambda \otimes {\mathcal L}_{\pi^{{\sharp}}\alpha}
\lambda \ ,
\end{equation}
for all 1-forms $\alpha$.
For any generator $\partial$ of the Gerstenhaber bracket $[~,~]_{\pi,
  \psi}$,
$$
[\alpha, \lambda]_{\pi, \psi} = (\partial \alpha) \lambda - \alpha
\wedge \partial \lambda \ .
$$
If $\partial = \partial _{\pi} + {\underline \partial}_{\pi, \psi}
+i_{Y_{\pi, \psi}}$, then
$\partial \alpha = <\alpha, Y_{\pi, \psi} > - i_{\pi} d \alpha$,
while,
by equation \eqref{*Z}, $\partial \lambda = -
i_{Z_{\pi, \psi,\lambda}} \lambda$,
and therefore
$$[\alpha, \lambda]_{\pi, \psi} =  (<\alpha,Z_{\pi, \psi,\lambda} +
Y_{\pi, \psi} > - i_{\pi} d \alpha) \lambda \ .$$
Since by definition $X_{\pi, \psi}$ satisfies
  ${\mathcal L}_{\pi^{{\sharp}}\alpha} \lambda = 
(<\alpha, X_{\pi, \psi}> + i_{\pi} d\alpha) \lambda $,
the vector field $U$ is such that
$$
<\alpha, U > \lambda \otimes \lambda = 2 <\alpha, Z_{\pi,\psi,
  \lambda} > \lambda \otimes \lambda \ .
$$
Therefore $U= 2 Z_{\pi,\psi,
  \lambda}$.
\end{proof}

In the Poisson case, $\partial = \partial_{\pi}$ and $\partial_{\pi}
\lambda = - i_{X_{\pi, \lambda}} \lambda$, and we obtain the relation
$
<\alpha, U > \lambda \otimes \lambda = 2 <\alpha, X_{\pi, \lambda} >
\lambda \otimes \lambda \ ,
$
which gives a new proof of the
result of \cite{ELW} stating that the modular class of a
Poisson manifold,
defined as the class of $X_{\pi,\lambda}$, characterized by an
equation such as \eqref{propX}, is equal to one half the ELW-modular class of
the Lie algebroid $(T^*M, \pi^{{\sharp}}, [~,~]_{\pi})$.

\medskip

The simple relationship between the modular class of $A=TM$
and the ELW-modular class of $A^*=T^*M$ does not in general hold for a Lie
algebroid with a twisted Poisson structure. In particular,
in Section \ref{section8}, we shall show that the two
classes may be different for a Lie algebra considered as a
Lie algebroid over a point, even in the usual, untwisted case.

\subsection{Modular class and gauge transformations}

Assume that we are given, on the cotangent bundle  $T^*M$
of a manifold $M$, two Lie algebroid structures,
denoted by $(T^*M, \rho, [\, , \,]) $ and $(T^*M,\rho', [\, , \,]') $
respectively, together with a Lie algebroid isomorphism,
$\sigma: T^*M \to T^*M$, over the identity of the base $M$.
According to formula  \eqref{representation}, each Lie algebroid
acts on $ L^{T^*M}= \wedge^n T^*M \otimes\wedge^n  T^*M$,
where $n$ is the dimension of $M$, and
 we denote
by $ D$ and  $D'$ these representations.
Define $\tau : L^{T^*M} \to L^{T^*M} $ by $\omega \mapsto \det(\sigma)
\, \omega$,
for all $ \omega \in \Gamma (L^{T^*M})$.

\begin{lemma} The isomorphism $  \tau$ is an intertwining operator
for the representations $ D$ and $ D' \circ \sigma$, i.e.,
\begin{equation}\label{eq:2D} D_{\sigma (\alpha)}' \tau (\omega)=\tau (D_\alpha \omega),  \end{equation}
for all $\alpha \in \Gamma(T^*M)$, $\omega \in \Gamma(L^{T^*M})$. \end{lemma}
\begin{proof}
For any sections $ \omega_1$ and $\omega_2$ of $ \wedge^n T^*M$,
$  \tau (\omega_1 \otimes \omega_2) = (\wedge^n \sigma) (\omega_1)
\otimes \omega_2 $. Hence
$$   D'_{\sigma (\alpha)} \tau (\omega_1 \otimes \omega_2) =
[\sigma (\alpha), (\wedge^n \sigma) (\omega_1)]' \otimes \omega_2
+ (\wedge^n \sigma) (\omega_1 ) \otimes {\mathcal L}_{ \rho' (\sigma
  (\alpha ))} \omega_2   $$
$$=(\wedge^n \sigma) ([\alpha, \omega_1]) \otimes \omega_2
+(\wedge^n \sigma) (\omega_1 ) \otimes {\mathcal L}_{ \rho (\alpha)} \omega_2
=\det(\sigma) D_\alpha (\omega_1 \otimes \omega_2) , $$
since $\sigma $ is an isomorphism of Lie algebroids.
\end{proof}

Choose a volume form $\lambda$ on $T^*M$.
Let $U$ (resp., $U'$) be the representative of the
ELW-modular class of
the Lie algebroid $(T^*M , \rho, [\, , \,])$ (resp.,
 $(T^*M, \rho', [\, , \,]')$) with volume form $\lambda$ (resp.
$\mu = \sqrt{|\det (\sigma)|} \lambda$).

\begin{proposition} \label{prop:isopullback} The vector
  fields $U$ and $U'$ are related by
$  U =   ^t \! \! \sigma \, \, U' $, where $^t \sigma$ is the
transpose of $\sigma$.
\end{proposition}
\begin{proof}
For some locally constant function
$\epsilon \in \{-1,+1\}$, $\tau (\lambda \otimes \lambda) = \epsilon \mu \otimes \mu   $.
Therefore, by equation \eqref{eq:2D},
$$  \epsilon D_{\sigma (\alpha)}' \mu \otimes \mu=\tau (D_\alpha \lambda \otimes \lambda)  .  $$
By the definition of the modular vector field given by equation \eqref{U},
$$  <  \sigma(\alpha),U' >  \epsilon \mu \otimes \mu =   < \alpha,U >
\tau (\lambda \otimes \lambda) . $$
Hence $ <  \sigma(\alpha),U' >=<  \alpha,U >$, and the result follows.
\end{proof}

A twisted Poisson structure
on a given manifold, $M$,
can be modified by a {\it gauge transformation}.
Assume that $B$ is a 2-form on $M$ such that
for all $m \in M$, the linear automorphism
of $T_m^*M$,  $\sigma_B : \alpha \mapsto \alpha + i_{\pi^{\sharp} \alpha} B$,
is invertible. Define a bivector $\pi'$ by
$$ (\pi')^{\sharp} =\pi^{\sharp}  ({\rm {Id}} + B^{\flat}  \circ \pi^{\sharp} )^{-1} \ , $$
where $B^{\flat}: TM \to T^*M$ is defined by $ B^{\flat}  (X)=i_X B$.
Then $(\pi',\psi-dB)$ is a twisted Poisson structure
that is said to be obtained by a {\em gauge transformation}
from $(\pi,\psi)$ \cite{severa}.
\begin{proposition}
A vector field $X$ on $M$ is in the modular class of $(TM,\pi,\psi) $
if and only if $X + \pi^{\sharp}  i_X B $
is in the modular class of $(TM,\pi',\psi-dB) $.
\end{proposition}
\begin{proof}
According to \cite{severa},
the map $\sigma_B$
is a Lie algebroid isomorphism
from $T^*M$ equipped with the Lie algebroid structure
associated to $(\pi,\psi)$
to $ T^*M$ equipped with the Lie algebroid structure
associated to $(\pi',\psi - dB) $.
 It follows from Propositions \ref{prop:evens} and
\ref{prop:isopullback} that its transpose, $X \mapsto X +
 \pi^{\sharp} i_X B $, maps modular
class to modular class.
\end{proof}

\section{Examples}\label{section8}

\noindent{\bf Example 1.
The case where $(\wedge^3 \pi^{{\sharp}})  \psi=0$.}
Whenever $ (M,\pi)$ is a Poisson manifold
and $\psi$ is a closed 3-form
satisfying $ (\wedge^3 \pi^{{\sharp}} ) \psi=0$,
then $(M,\pi,\psi) $ is a twisted Poisson manifold.
One can prove that the assumption $(\wedge^3 \pi^{{\sharp}})  \psi=0$
implies $Y_{\pi,\psi}=0$.
Hence, for any volume form $\lambda$,
the modular vector field of $(M,\pi,\psi)$ is
equal to $X_{\pi,\lambda}$, which is the modular
vector field of the Poisson manifold $(M,\pi)$.
In this case therefore,
the modular class of the twisted Poisson manifold
$ (M,\pi,\psi)$ is equal to the modular class
of the Poisson manifold $(M,\pi)$.
This conclusion holds more generally in the case of a Lie algebroid
with a twisted Poisson structure of this type.

\medskip

\noindent{\bf Example 2. The case where $\pi ^{{\sharp}}$ is invertible.}
Let $(\pi,\psi)$ be a twisted Poisson structure on $A$ such that $\pi \in \Gamma (\wedge^2 A)$
is of maximal rank at each point. Let $\omega$ be the
non-degenerate
2-form such that $\omega^{\flat} = (\pi^{{\sharp}})^{-1}$.
One can show that the modular vector field $Z_{\pi,\psi,\lambda}$, where
$\lambda$ is the volume form $\omega^{\frac{N}{2}}$,
vanishes. The proof rests on the properties $d_A \omega = \psi$ and
$\frac{N}{2} \omega^{\frac{N}{2}-1}=
 i_{\pi}\omega^{\frac{N}{2}}$.
Therefore the Lie algebroids with a
twisted Poisson structure whose bivector is of maximal rank are
unimodular, a result which extends the fact that symplectic
manifolds are unimodular.

\medskip

\noindent{\bf Example 3. Twisted Poisson structures on Lie groups.}
Consider the example of a
twisted Poisson structure defined on a dense open
subset of a Lie group  \cite{severa}.
Let $G$ be a Lie group with Lie algebra
${\mathfrak g}$, and assume that $\mathfrak g$ is
equipped with an invariant non-degenerate
symmetric bilinear form $\big< \ , \ \big>$.
For any $X \in {\mathfrak g}$ (resp., $ \alpha \in {\mathfrak g}^*$),
we denote by $ \overline{X} \in {\mathfrak g}^*$
(resp., $  \underline{\alpha} \in {\mathfrak g}$) the image of $X$
(resp.,~$\alpha$) under the isomorphism of
$ {\mathfrak g}$ to ${\mathfrak g}^*$
(resp., $ {\mathfrak g}^*$ to ${\mathfrak g}$)
induced by $\big< \ , \ \big>$.
For any $ X \in \wedge^\bullet {\mathfrak g}$,
let ${X}^R$ and ${X}^L$
be the corresponding right- and left-invariant
multivector fields.
Dually, for any $ \alpha \in \wedge^{\bullet} {\mathfrak g}^*$,
we denote by $ \alpha^R$ and $\alpha^L$
the corresponding right- and left-invariant forms.
Then, for any $\alpha \in {\mathfrak g}^*$,
$ d \alpha^L = ( d_{\mathfrak g} \alpha )^L$,
while $ d \alpha^R = -( d_{\mathfrak g} \alpha )^R$,
where $d_{\mathfrak g}$
is the Chevalley-Eilenberg differential of $\mathfrak g$
and $d$ is the de Rham differential.

The canonical 3-form $\psi$ on $\mathfrak g$, defined by
$ \psi (X,Y,Z) = \frac{1}{2} \big< X , [ Y, Z ]  \big>$,
for $X, Y$ and $Z \in \mathfrak g$,
satisfies $ \psi^R = \psi^L$. The corresponding bi-invariant form
on $G$ is called the Cartan 3-form, and we shall denote it
by the same symbol. It satisfies, for any $X \in \mathfrak g$,
\begin{equation} \label{eq:contpsi}
i_{{X}^L} \psi  = - \frac{1}{2}
 (  d_{\mathfrak g} \overline{X} )^L  \quad \quad {\mathrm
  {and}}  \quad \quad i_{ {X}^R} \psi   = -\frac{1}{2}
  (  d_{\mathfrak g} \overline{X})^R .
\end{equation}

Let $G_0$
be the dense open set of elements $g \in G$
such that $ -1$ is not an eigenvalue
of ${\rm {Ad}}_g$.
Equivalently, $G_0$
is the subset of elements $g \in G$
such that the linear map
$\alpha \mapsto \alpha^R(g) + \alpha^L (g)$
is an isomorphism from $ {\mathfrak g}^*$ to $T^*_g G $.
A twisted Poisson structure on $G_0$
is given by $\psi$
and the bivector $\pi$ defined by
$$
 \pi^{{\sharp}} (\alpha^L + \alpha^R)= 2 ({\underline {\alpha}}^L
- {\underline {\alpha}}^R)
\ ,
$$
\begin{proposition}
The twisted Poisson manifold $(G_0,\pi,\psi)$ is unimodular.
\end{proposition}
\begin{proof}
We first evaluate
$Y_{\pi,\psi} $
on the 1-forms
$\alpha^R + \alpha^L$, obtaining
$$ i_{Y_{\pi,\psi}}(\alpha^L + \alpha^R)=
  -  i_{ \pi \wedge \pi^{{\sharp}} (\alpha^L + \alpha^R)  } \psi
=  - 2   i_{\pi \wedge ({\underline {\alpha}}^L
- {\underline {\alpha}}^R)} \psi .
$$
By equation (\ref{eq:contpsi}), we obtain
\begin{equation}  \label{eq:step1caseG} i_{Y_{\pi,\psi}}(\alpha^L +
  \alpha^R) = i_\pi ((d_{\mathfrak g} \alpha)^L - (d_{\mathfrak g}\alpha)^R ) .
\end{equation}
We can choose a volume
form $\lambda$ which is  both left- and right-invariant.
By definition,
$$ (i_{X_{\pi,\lambda}} (\alpha^L + \alpha^R)) \, \lambda
= {\mathcal L}_{ \pi^{{\sharp}} (\alpha^L + \alpha^R) } \lambda -
(i_\pi  d  (\alpha^L + \alpha^R) ) \, \lambda \ .
$$
Since $\lambda$ is left- and right-invariant,
${\mathcal L}_{ \pi^{{\sharp}} (\alpha^L + \alpha^R) } \lambda = 
2 ( {\mathcal L}_{{\underline{\alpha}}^L} \lambda -
{\mathcal L}_{{\underline{\alpha}}^R} \lambda ) = 0$, thus
\begin{equation}\label{eq:step3caseG}   i_{X_{\pi,\lambda}} (\alpha^L + \alpha^R)=  - i_\pi ( ( d_{\mathfrak g}\alpha)^L - (d_{\mathfrak g}\alpha)^R ) .\end{equation}
Equations (\ref{eq:step1caseG}) and (\ref{eq:step3caseG}) imply that
$Z_{\pi,\psi,\lambda}=0 $.
\end{proof}

\medskip

\noindent{\bf Example 4. Lie algebras with a triangular $r$-matrix.}
Let  $\mathfrak g$ be a real Lie algebra of dimension $N$,
and let $r\in \wedge^2 {\mathfrak g}$ such that $[r,r] =0$, this bracket
being the algebraic Schouten bracket on $\wedge^{\bullet} \mathfrak
g$. Such an $r$ is called a {\em triangular $r$-matrix}.
Then ${\mathfrak g}^*$ is a Lie algebra with bracket $[\alpha,
\beta]_r = 
{\rm {ad}}^*_{ r^{{\sharp}} \alpha } \beta -{\rm {ad}}^*_{r^{{\sharp}} \beta } \alpha$.
Let $\lambda$ be a non-zero element in $\wedge^{N}{\mathfrak g}^*$.
The class that we have just defined is the class of the
element $X_{r, \lambda}$ of $\mathfrak g$ such that
$$
  <\alpha,  X_{r, \lambda}> \lambda = (\partial_r   \alpha) \lambda +d_{\mathfrak g} (i_{r^{{\sharp}}\alpha} \lambda)  ,
$$
or, equivalently, by equation \eqref{*X},
$$  <\alpha,  X_{r, \lambda}> \lambda = - \alpha \wedge \partial_r \lambda  ,
$$
for $\alpha \in {\mathfrak g}^*$. On the other hand,
the ELW-modular class of ${\mathfrak g}^*$ considered as a Lie
algebroid over a point is the class,
in the Lie algebra cohomology of $\mathfrak
g$, of the element $\widetilde X_{r, \lambda}$ of $\mathfrak g$
such that
$$
<\alpha, \widetilde X_{r, \lambda}> \lambda = [\alpha, \lambda]_r .
$$
Since $\partial_r$ generates the bracket $[~,~]_r$,
$[\alpha, \lambda]_r = (\partial_r \alpha) \lambda - \alpha \wedge \partial_r
\lambda$. Thus the class of $X_{r,\lambda}$ is one half that of $\widetilde
X_{r, \lambda}$
if and only if
$d_{\mathfrak g} (i_{r^{{\sharp}}\alpha} \lambda)$ vanishes for all $\alpha
\in {\mathfrak g}^*$.

\medskip

\noindent{\it Example 4.1. Triangular $r$-matrix on
the non-abelian 2-dimensional Lie algebra.}
In this simple example, one class vanishes, while the other does not.
Let $\mathfrak g$ be the non-abelian 2-dimensional Lie algebra
with basis $(e_1, e_2)$ such that $[e_1,e_2]= e_1$,
with the triangular $r$-matrix $\pi= r= e_1 \wedge e_2$, and $\psi=0$.
For the dual basis
$(\epsilon^1, \epsilon^2)$, $[\epsilon^1, \epsilon^2]_r = -
\epsilon^2$, the dual Lie algebra is in fact isomorphic to $\mathfrak g$.
Let $\lambda = \epsilon^1 \wedge \epsilon^2$. The
modular vector $X_{r, \lambda}$ vanishes
  since $i_{X_{r,\lambda}}\lambda  = -  d_{\mathfrak g}i _{r} \lambda$
  and $d_{\mathfrak g}$
vanishes on $\wedge^0{\mathfrak g}$.
On the other hand, the ELW-modular class of $\mathfrak g$ is the class
of
$
{\widetilde X}_{r ,\lambda}= - e_1$, and this element is non-zero,
therefore its class does not vanish.

\medskip

\noindent{\it Remark.}
The ELW-modular class of ${\mathfrak g}^*$
is equal to the class of the infinitesimal modular
character of ${\mathfrak g}^*$, which is the element of ${\mathfrak g}$,
$\alpha \mapsto {\rm {Tr}}({\rm {ad}}_\alpha)$, and
it also coincides with the modular class of the linear
Poisson manifold $\mathfrak g$, dual to the Lie algebra ${\mathfrak
  g}^*$.  In the example above, alternatively, we can
show that the constant vector field ${\widetilde X}_{r
  ,\lambda}= - e_1$
is not globally
hamiltonian with respect to the linear Poisson structure $\pi^r$
on $\mathfrak g$, the dual of the Lie algebra $({\mathfrak g}^*,
[~,~]_r)$. If $(x_1,x_2)$
are the coordinates on $\mathfrak g$ with respect to the chosen
basis,
and $x \in {\mathfrak g}$, then
$\pi^r(x) = - {x_2} e_1 \wedge e_2 $. Since, for a function
$u \in C^{\infty}({\mathfrak g})$,
$d_{\pi^r}u
 = [\pi^r,u] = {x_2}(\frac{\partial u}{\partial x_2}
  e_1 - \frac{\partial u}{\partial x_1} e_2)$, the condition
  $-e_1 = d_{\pi^r}u$ is clearly not realizable along the axis
  $x_2 =0$.

\medskip

\noindent{\it Example 4.2. Triangular $r$-matrix on
${\mathfrak {sl}}(2,\mathbb R)$.}
On ${\mathfrak g}={\mathfrak {sl}}(2,{\mathbb R})$,
define $\pi=r =X_+ \wedge H $,
where $X_+$, $X_-$ and $H$ denote the usual basis of $ {\mathfrak
 {sl}}_2({\mathbb R})$ such that
$[H,X_+]=2X_+,\,  [H,X_-] = -2 X_-$ and $[X_+,X_-]=H $.
The relation $[r,r]=0$ is satisfied, hence
$(r,\psi )$ with $\psi =0$ is a twisted Poisson structure on
$ {\mathfrak {sl}}(2,{\mathbb R})$. In view of equation \eqref{*Z},
in the particular case where $\psi =0$,
for any non-zero $\lambda \in \wedge^3 {\mathfrak g}^*$,
the modular vector ${Z_{r, \lambda}}$
is  defined by  $*_\lambda {Z_{r, \lambda}}
=  - d_{{\mathfrak g}} i_r \lambda  $,
We choose $ \lambda = H^* \wedge X_+^* \wedge X_-^* $,
where $X_+^*,X_-^* , H^*$ is the dual basis, and we compute
$$  * _\lambda Z_{r, \lambda} = - d_{{\mathfrak g}} i_r
  \lambda = -
 d_{ {\mathfrak g} } X_-^* =  - 2 H^* \wedge  X_-^*= i_{2 X_+} \lambda .$$
Hence the modular vector is $Z_{r,\lambda}=2 X_+$
and the modular class is not trivial.

It is straightforward to check that
$d_{{\mathfrak g}}  i_{r^{\sharp} \alpha} \lambda$
vanishes for any $\alpha \in   {\mathfrak {sl}}(2,{\mathbb R})^* $.
Therefore, according to the above statement concerning the
case of triangular $r$-matrices,
the ELW-modular class of the Lie algebra ${\mathfrak g}^*$ is equal
to twice the modular class. This conclusion can be verified by
computing the infinitesimal modular character of the Lie algebra ${\mathfrak
  g}^*$.

\medskip

\noindent{\bf Example 5. A twisted $r$-matrix.}
Let ${\mathfrak g} $ be the Lie algebra of the group of
affine transformations of
${\mathbb R}^2$. Denote by $u_1,u_2$
a basis of ${\mathbb R}^2$ (that we identify with
the abelian subalgebra of translations) and
$e_{i,j}, \, i,j \in \{1,2\} $,
 the  basis of $ {\mathfrak {gl}}(2,{\mathbb R})$ given by
$e_{i,j} ( u_k ) =[e_{i,j},u_k]= \delta_{jk} u_i $.
The dual basis is denoted by $(e_{1,1}^*, e_{1,2}^*,
e_{2,1}^*,e_{2,2}^*, u_1^*, u_2^*)$.

Define $\pi =r =  e_{1,1} \wedge e_{2,2} + u_1 \wedge u_2 $
and $ \psi =
-(e_{1,1}^* + e_{2,2}^* ) \wedge u_1^* \wedge u_2^* $.
It is easy to check that $\psi$
is closed, and the following computation shows that
$(r,\psi)$ is a twisted Poisson structure on $\mathfrak g$,
$$ \frac{1}{2}[r,r]= [e_{1,1} \wedge e_{2,2},u_1 \wedge u_2  ]=(e_{1,1}-
e_{2,2} )\wedge u_1 \wedge u_2 =-(\wedge^3 r^{\sharp})( (e_{1,1}^* + e_{2,2}^*
)  \wedge u_1^* \wedge u_2^*).$$
To compute the modular field, we first evaluate
$$ Y_{\pi,\psi}  =   r^{\sharp} (e_{1,1}^* + e_{2,2}^*) = e_{2,2} - e_{1,1}   .$$
Then we set
$ \lambda = e_{2,2}^* \wedge  e_{1,1}^* \wedge e_{1,2}^*
\wedge e_{2,1}^* \wedge u_2^* \wedge u_1^* $,
and we obtain
$$ X_{\pi,\lambda} = -*_{\lambda}^{-1} d_{\mathfrak g}
i_{r} \lambda =   -*_{\lambda}^{-1} d_{\mathfrak g}  (   e_{1,2}^*
\wedge e_{2,1}^* \wedge u_2^* \wedge u_1^* +    e_{2,2}^* \wedge  e_{1,1}^* \wedge e_{1,2}^*
\wedge e_{2,1}^*  )   $$
$$ = - *_{\lambda}^{-1} (e_{1,2}^* \wedge e_{2,1}^* \wedge
d_{\mathfrak g} (u_2^* \wedge u_1^*)   )  =   e_{2,2} - e_{1,1}  .
$$
Hence the modular class of $({\mathfrak g}, r, \psi)$ is
$2(e_{2,2} - e_{1,1})$.

\medskip

\noindent{\bf Remark.}
Whereas to a Lie algebroid with a Poisson structure, $(A,\pi)$,
is associated a Lie
bialgebroid structure on $(A,A^*)$, to a Lie algebroid
with a twisted Poisson
structure, $(A, \pi, \psi)$, is associated a {\em quasi-Lie
bialgebroid} structure on $(A,A^*)$. This was proved by Roytenberg in \cite{R}.
The Lie algebroid structure on $A^*$ is that described in Theorem
\ref{dualalgd}, while the bracket on $A$ is $[~,~]_A + \psi^{(1)}$,
and the associated derivation of $\Gamma (\wedge^{\bullet}A^*)$ is the
operator $d_{A,\pi,\psi} =d_A + \delta_{\pi,\psi}$ introduced in
\eqref{quasidiff}.
If, in particular, $A$ is a Lie algebra $\mathfrak g$, and $r \in \wedge
^2 \mathfrak g$ and a $d_{\mathfrak g}$-cocycle 
$\psi \in \wedge
^3 {\mathfrak g}^*$ satisfy 
the twisted Poisson condition, $\frac{1}{2}[r,r] =
(\wedge^3 r^{{\sharp}})\psi$, 
then the pair $({\mathfrak g},{\mathfrak g^*})$
becomes a {\em quasi-Lie bialgebra}.
To avoid confusion, we stress that
the term ``twisted'' is used here to denote an
$r$-matrix which is not triangular but satisfies a non-linear condition
involving a 3-cocycle, while in the theory of Lie
quasi-bialgebras,
a``twist'' or ``twisting'' is the modification of an $r$-matrix
by the addition of an element $t \in \wedge^{2}{\mathfrak g}$.
The case of ``twisted $r$-matrices'' in the former sense deserves to
be further explored.

\vspace{0.5cm}

\noindent{\bf Acknowledgments}
We are grateful to Alan Weinstein for his very valuable comment on the
non-orientable case, and we thank the referee for his remarks.

\bigskip

\bigskip

\bigskip

\noindent
      Yvette Kosmann-Schwarzbach\\
      Centre de Math\'ematiques Laurent Schwartz, U.M.R. du CNRS 7640\\
      \'Ecole Polytechnique\\
      F-91128 Palaiseau\\
      yks@math.polytechnique.fr

\medskip

\noindent
Camille Laurent-Gengoux \\
D{\'e}partement de Math{\'e}matiques\\
Universit{\'e} de Poitiers \\
F-86962 Futuroscope Cedex\\
Camille.Laurent@math.univ-poitiers.fr

\label{lastpage}
\end{document}